\documentclass[12pt]{article}
\usepackage{amsmath,amsfonts,amssymb}
\usepackage[all]{xy}
\usepackage{graphicx}
\usepackage{epsf}
\usepackage{cancel}
\usepackage{ulem}
\usepackage{multirow}
\usepackage{rotating}
\usepackage{arydshln}
\usepackage{tikz}
\usepackage{eucal}
\usepackage[colorlinks=true,urlcolor=black]{hyperref}
\usetikzlibrary{decorations.markings}
\usetikzlibrary{arrows}
\usepackage[colorlinks=true,urlcolor=black]{hyperref}
\hypersetup{			
backref = true,			
pagebackref = true,		
hyperindex = true, 		
colorlinks = true, 		
breaklinks = true, 		
urlcolor = blue, 		
linkcolor = blue, 		
bookmarks = true,		
bookmarksopen = true,  
citecolor=red,
}

\advance \textheight by \topskip
 

\newtheorem{theorem}{Theorem}[section]

\newtheorem{corollary}[theorem]{Corollary}

\newtheorem{definition}[theorem]{Definition}
\newtheorem{example}[theorem]{Example}

\newtheorem{lemma}[theorem]{Lemma}

\newtheorem{proposition}[theorem]{Proposition}
\newtheorem{remark}[theorem]{Remark}

\newtheorem{liste}[theorem]{List}

\newcommand{\qeed}{\hfill\textrm{QED}\break\null}
\newenvironment{demo}{\noindent\textit{Proof.}~}{\qeed}

\def\F{{\cal F}}

\def\FF{\mathbb F}

\newcommand{\beq}{\begin{equation}}
\newcommand{\eeq}{\end{equation}}
\newcommand{\beqa}{\begin{eqnarray}}
\newcommand{\eeqa}{\end{eqnarray}}

\newcommand{\noi}{\noindent}
\newcommand{\nn}{\nonumber}

\newcommand{\e}{\epsilon}
\def\>{\rangle}
\def\<{\langle}

\begin{document}

\title{Incidence geometry of the Fano plane and Freudenthal's ansatz for the  construction of     (split) octonions.}
{\bf }

\author{
{\sf M. Rausch de Traubenberg}\thanks{e-mail:
michel.rausch@iphc.cnr.fr}$\,\,$${}^{a}$ and
{\sf M. J. Slupinski}\thanks{e-mail:
marcus.slupinski@math.unistra.fr}$\,\,$${}^{b}$
\\
{\small ${}^{a}${\it IPHC-DRS, UdS, CNRS, IN2P3 
}}\\
{\small {\it  23  rue du Loess, Strasbourg, 67037 Cedex, France 
}}\\
{\small ${}^{b}${\it
Institut de Recherches en Math\'emathique Avanc\'ee,
UdS and CNRS}}\\
{\small {\it 7 rue R. Descartes, 67084 Strasbourg Cedex, France.}}  \\ 
}

\maketitle
\date

\begin{abstract}
\noindent{\small
In this article we consider  structures on a Fano plane  $\F$  which allow a generalisation of Freudenthal's construction
of a {norm and a bilinear }multiplication law  on an eight-dimensional vector  space $\mathbb O_\F$ canonically associated to $\F$.
{We first determine
necessary and sufficient   conditions  in terms of the incidence geometry of $\F$ for these structures to give rise to  division composition
algebras, and classify   the corresponding structures
using a  logarithmic version  of the multiplication}.  We then show how these results can be used to deduce analogous results in the { split composition} algebra case.

}
\end{abstract}

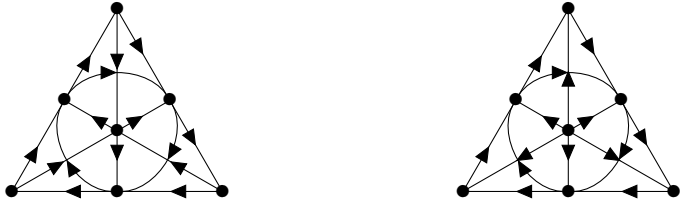
\begin{figure}[!ht]
  \begin{center}\hskip 2.truecm
  \begin{minipage}{4cm}
    \begin{tikzpicture}[scale=.7]
\tikzstyle{point}=[circle,draw]

\tikzstyle{ligne}=[thick]
\tikzstyle{pointille}=[thick,dotted]

\tikzset{->-/.style={decoration={
  markings,
  mark=at position .5 with {\arrow{>}}},postaction={decorate}}}

\tikzset{middlearrow/.style={
        decoration={markings,
            mark= at position 0.5 with {\arrow{#1}} ,
        },
        postaction={decorate}
    }}

\coordinate  (3) at ( -2,-1.15);

\coordinate  (2) at ( 0,-1.15) ;
\coordinate  (5) at ( 2,-1.15);

\coordinate  (6) at (0,2.31);
\coordinate  (1) at (1, .58);

\coordinate  (4) at (-1., .58) ;
\coordinate  (7) at (0,0) ;

\draw [color=black, middlearrow={triangle 45}] (1) to [bend left=65] (2);
\draw [color=black,  middlearrow={triangle 45}] (2) --(3);
\draw [color=black, middlearrow={triangle 45}] (3) --(4);
\draw [color=black,  middlearrow={triangle 45}] (6) --(7);
\draw [color=black,  middlearrow={triangle 45}] (7) --(1);

\draw [middlearrow={triangle 45}] (6) --(1);
\draw [ middlearrow={ triangle 45}] (1) --(5);

\draw[middlearrow={triangle 45}]  (5) --(2);

\draw[middlearrow={triangle 45}]  (3) --(7);

\draw  [color =black, middlearrow={triangle 45}]  (5) --(7);
\draw [color =black, middlearrow={triangle 45}] (7) --(4);
\draw[middlearrow={triangle 45}]  (4) --(6);
\draw [middlearrow={triangle 45}]  (7) --(2);

\draw [middlearrow={triangle 45}] (2) to [bend left=65] (4);
\draw [middlearrow={triangle 45}]  (4) to [bend left=65] (1);
\draw (1)  node[scale=0.6] {$\bullet$} ;
\draw (1)  node {$\bullet$} ;

\draw (2)  node{$\bullet$} ;
\draw (3) [color=black]  node{$\bullet$} ;
\draw (4)  node{$\bullet$} ;
\draw (5)  node{$\bullet$} ;
\draw (6)  node{$\bullet$} ;
\draw (7)  node{$\bullet$} ;

  \end{tikzpicture} 
  \end{minipage}\hskip 2.truecm
  \begin{minipage}{8cm}
 \begin{tikzpicture}[scale=.7]
\tikzstyle{point}=[circle,draw]

\tikzstyle{ligne}=[thick]
\tikzstyle{pointille}=[thick,dotted]

\tikzset{->-/.style={decoration={
  markings,
  mark=at position .5 with {\arrow{>}}},postaction={decorate}}}

\tikzset{middlearrow/.style={
        decoration={markings,
            mark= at position 0.5 with {\arrow{#1}} ,
        },
        postaction={decorate}
    }}

\coordinate  (3) at ( -2,-1.15);

\coordinate  (2) at ( 0,-1.15) ;
\coordinate  (5) at ( 2,-1.15);

\coordinate  (6) at (0,2.31);
\coordinate  (1) at (1, .58);

\coordinate  (4) at (-1., .58) ;
\coordinate  (7) at (0,0) ;

\draw [color=black, middlearrow={triangle 45}] (1) to [bend left=65] (2);
\draw [color=black,  middlearrow={triangle 45}] (2) --(3);
\draw [color=black, middlearrow={triangle 45}] (3) --(4);
\draw [color=black,  middlearrow={triangle 45}] (7) --(6);
\draw [color=black,  middlearrow={triangle 45}] (7) --(1);

\draw [middlearrow={triangle 45}] (6) --(1);
\draw [ middlearrow={ triangle 45}] (1) --(5);

\draw[middlearrow={triangle 45}]  (5) --(2);

\draw[middlearrow={triangle 45}]  (7) --(3);

\draw  [color =black, middlearrow={triangle 45}]  (7) --(5);
\draw [color =black, middlearrow={triangle 45}] (7) --(4);
\draw[middlearrow={triangle 45}]  (4) --(6);
\draw [middlearrow={triangle 45}]  (7) --(2);

\draw [color=black, middlearrow={triangle 45}] (2) to [bend left=65] (4);
\draw [color=black,middlearrow={triangle 45}]  (4) to [bend left=65] (1);
\draw (1)  node[scale=0.6] {$\bullet$} ;
\draw (1)  node[color=black] {$\bullet$} ;

\draw (2)  node[color=black]{$\bullet$} ;
\draw (3) [color=black]  node{$\bullet$} ;
\draw (4)  node[color=black]{$\bullet$} ;
\draw (5)  node{$\bullet$} ;
\draw (6)  node{$\bullet$} ;
\draw (7)  node{$\bullet$} ;

  \end{tikzpicture}
\end{minipage}
  
  \caption{Systems of arrows giving division and split composition algebras.}
\label{fig:oct}
\end{center}
\end{figure}

{In  \cite{freu} (p.19) } Freudenthal  gave an  {ansatz to construct } a real  eight-dimensional   division composition algebra from
  {a particular  system of arrows (Figure
\ref{fig:oct} left panel) between pairs of points in  a Fano plane $\F$} (see also
\cite{vdb,zorn,SV,baez} for a historical perspective).
 { A slight variant of this ansatz, whose origin in the literature is uncertain,  allows one to construct  eight-dimensional  split composition algebras from other  systems of arrows in $\F$ (for example, Figure
\ref{fig:oct} right panel).} { The main aim of  this short paper is to characterise in terms of the geometry of $\F$  those systems of arrows in $\F$ which give rise to composition algebras by   Freudenthal's ansatz or its variant.}


{We begin by defining the notions of  norm and multiplication factor  (a formal version
of a system of arrows) on $\F$,   and from them construct a quadratic form and bilinear multiplication law  on an eight-dimensional vector  space $\mathbb O_\F$ canonically associated to $\F$. 
This construction can  give a division composition algebra if and only if the norm is constant,
and in that case we show (see Theorems  \ref{theo:future}  and \ref{theo:future1}) that necessary and sufficient conditions for obtaining  a  division composition algebra   are that
 either the future of every point is a  line or the past of every point is a line (see Definition \ref{def:norm}
 for the past/future of points with respect to a multiplication factor).}
By taking {logarithms }  (see Definition \ref{def:exp})
we classify all such multiplication factors {and it turns out  there are eight of each type}.
From a different geometric point of view 
we  show that any  set of seven triangles in $\F$ satisfying the axioms of a Fano plane determines a unique
division composition algebra such that the given set of triangles is exactly the set of pasts of points (see Theorem \ref{theo:future2}).

{If a norm is not constant  it defines a line in $\F$  and we use this observation to establish a connection between the split composition algebra case and the division  composition algebra case of the above procedure.  More precisely, we show there is  a  bijection between the set of multiplication factors giving rise to  split composition algebras with the same norm and the set of multiplication factors giving rise to division composition algebras.} This allows us to characterise all multiplication factors which give rise to split composition algebras.

\section{The Fano plane: definition and  basic properties}
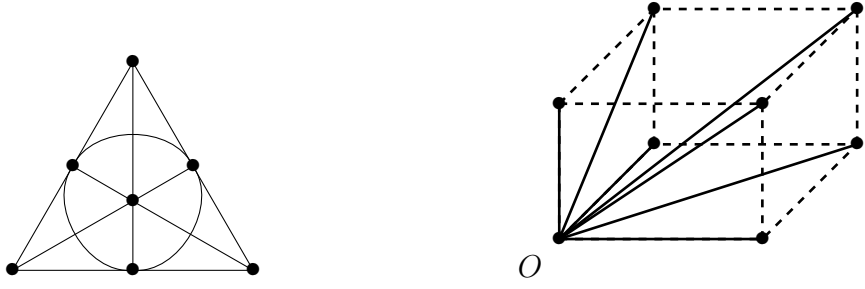
\begin{figure}[!ht]
  \begin{center}

\vskip -1truecm
    
  \begin{tikzpicture}[scale=.8] 
\tikzstyle{point}=[circle,draw]

\tikzstyle{ligne}=[thick]
\tikzstyle{pointille}=[thick,dotted]

\tikzset{->-/.style={decoration={
  markings,
  mark=at position .5 with {\arrow{>}}},postaction={decorate}}}

\tikzset{middlearrow/.style={
        decoration={markings,
            mark= at position 0.5 with {\arrow{#1}} ,
        },
        postaction={decorate}
    }}

\coordinate  (3) at ( -2,-1.15);

\coordinate  (2) at  ( 0,-1.15) ;
\coordinate  (5) at   ( 2,-1.15);

\coordinate  (6) at   (0,2.31);
\coordinate  (1) at  (1, .58);

\coordinate  (4) at  (-1., .58) ;
\coordinate  (7) at  (0,0) ;

\draw  (6) --(1);
\draw (1) --(5);

\draw  (5) --(2);
\draw  (2) --(3);

\draw  (3) --(4);
\draw  (4) --(6);

\draw  (6) --(7);
\draw  (7) --(2);

\draw (3) --(7);
\draw  (7) --(1);

\draw (5) --(7);
\draw  (7) --(4);

\draw (1) to [bend left=65] (2);
\draw  (2) to [bend left=65] (4);
\draw  (4) to [bend left=65] (1);
\draw (1)  node {$\bullet$} ;

\draw (2)  node{$\bullet$} ;

\draw (3)  node{$\bullet$} ;
\draw (4)  node{$\bullet$} ;

\draw (5)  node{$\bullet$} ;

\draw (6)  node{$\bullet$} ;

\draw (7)  node{$\bullet$} ;

  \end{tikzpicture} \hskip 3.truecm 
    \begin{tikzpicture}[thick,scale=.6]
  \coordinate (P0) at (0,0);
    \coordinate (P1) at (4.5,0) ; 
    \coordinate (P4) at (6.6,2.1); 
    \coordinate (P2) at (2.1,2.1);

     \coordinate (P3) at (0,3);
    \coordinate (P7) at  (4.5,3) ;
    \coordinate (P6) at (6.6,5.1); 
    \coordinate (P5) at (2.1,5.1);

 \draw[solid][line width=1pt] (P0) -- (P1);
 \draw[solid][line width=1pt] (P0) -- (P2);
  \draw[solid][line width=1pt] (P0) -- (P3);
  \draw[solid][line width=1pt] (P0) -- (P4);
   \draw[solid][line width=1pt] (P0) -- (P5);
    \draw[solid][line width=1pt] (P0) -- (P6);
 
      \draw[solid][line width=1pt] (P0) -- (P7);

\draw[dashed][line width=1pt] (P0) -- (P1)-- (P4) -- (P2)--(P0);
    
    \draw[dashed][line width=1pt] (P3) -- (P7)-- (P6) -- (P5)--(P3);
    \draw[dashed][line width=1pt] (P0) --(P3);
    \draw[dashed][line width=1pt] (P1) --(P7);
    \draw[dashed][line width=1pt] (P4) --(P6);
    \draw[dashed][line width=1pt] (P2) --(P5);

\draw (P0)  node {$\bullet$} ;
\draw (P1)  node {$\bullet$} ;

\draw (P2)  node{$\bullet$} ;

\draw (P3)  node{$\bullet$} ;
\draw (P4)  node{$\bullet$} ;

\draw (P5)  node{$\bullet$} ;

\draw (P6)  node{$\bullet$} ;

\draw (P7)  node{$\bullet$} ;

   \draw (P0) node[below left=2pt]{$O$} ;

   \end{tikzpicture}
  
\caption{The Fano plane ${\cal F}$ and the Fano cube $V_{\cal F}$.}
\label{fig:FP}
 \end{center}
\end{figure}

\begin{definition}\label{def: Fano plane}
  \begin{enumerate}
\item A Fano plane is a set $\F$ with seven elements, called points, together with a   set $\F^\ast$ of seven
subsets of $\F$ containing  three elements, called lines. Points and lines are required to satisfy the axioms of a projective plane:
\begin{enumerate}
\item  two distinct points are contained in a unique line;
\item two distinct lines intersect in a unique point.
\end{enumerate}
It is a consequence of these properties that each point is contained in exactly three lines and three distinct lines intersecting in a point are called concurrent.
\item  The Fano cube $V_\F$ associated to $\F$ is the set $V_\F = \F \cup \{0\}$ with the unique $\mathbb Z_2-$vector space structure
such that if $P\ne Q$:
\beqa
P+Q=R \ \mbox{where} \ \{P,Q,R\} \ \mbox{is the unique line containing}\  P\  \mbox{and} \ Q\ .\nn
\eeqa

\item The Fano plane dual to $\F$ is the set  $\F^\ast$ whose ``points'' are the  lines of $\F$ and whose ``lines''
are triples of  concurrent lines in $\F$.

\end{enumerate}
\end{definition}
The Fano plane and the Fano cube are usually represented as in Figure \ref{fig:FP}.
\vskip .4truecm

\noi
 If $L=\{P,Q,R\}$ is triple of $\F$ then $L$ is a line {\it iff} $P+Q+R=0$.
A triangle is a subset   consisting of three points of $\F$ that  is not a line. 
A quadrilateral  is a subset consisting of four points of $\F$ of which no three are colinear. Equivalently a quadrilateral 
is the complement of a line.

\begin{definition}
Let $\Delta =\{P, Q, R\}$ be a triangle in $\F$. The orthopoint of $\Delta$ is $P+Q+R$ and the ortholine of $\Delta$ is
$\{P+Q, Q+R,R+P\}$. 
\end{definition}

We will need the following incidence properties:
\begin{liste} \label{li:inci}
\begin{enumerate}
\item[L1] If $\Delta = \{P,Q,R\}$ is a triangle then $\{P,Q,R,P+Q+R\}$ is a quadrilateral.
\item[L2] If $\Delta$ is a triangle,  $P_\Delta$ its orthpoint and $L_\Delta$ its ortholine, then $\F=P_\Delta\cup L_\Delta\cup \Delta$ is a disjoint union.
\item[L3] If $\F$ is a disjoint union of a point $P$, a line $L$  and a triangle $\Delta$,  then $P$ is the orthopoint of $\Delta$ and $L$ is the ortholine of $\Delta$.
\item[L4] Two distinct quadrilaterals intersect in two points.
\item[L5] The union of three distinct lines is either $\F$ if the lines are concurrent or the complement of a point if not.
\item[L6] Given two distinct points there are exactly two lines which contain neither.
\item[L7] If $\{P,Q,R,S\}$ is a quadrilateral then $P+Q+R+S=0$.
\end{enumerate}
\end{liste}

\section{Octonions and split octonions}
Let  $\FF$ be a field of characteristic not two and let $\F$ be a Fano plane equipped with a norm  and multiplication factor (see below). 
In this section, following Freudenthal \cite{freu}, we define a quadratic form  and a bilinear multiplication  on an  eight-dimensional $\FF-$vector space  ${\mathbb O_\F}$ canonically associated to $\F$. 
We give necessary and sufficient conditions for this procedure to define a  composition algebra, and prove that in the  case of a division composition algebra
either the future of every point is a line or the past of every point is a line (Theorem \ref{theo:future}).


\begin{definition}
\begin{enumerate}\label{def:norm}
\item A norm on $\F$ is a function $N: \F \to \{-1,1\}$  such that
\beqa
N(P+ Q) = N(P) N(Q) \ , \ \ P \ne Q \ .\nn
\eeqa
\item Let $\F^2_0 = \{(P,Q) \in \F^2 \ \ \text{s.t} \ \ P \ne Q \}$. 
A multiplication factor is a map $\epsilon: \F^2_0 \to \{-1,1\}$  such that
$\epsilon_{PQ}+\epsilon_{QP}=0$. For $P \in \F$ the future  (resp. past) $\overrightarrow{P_\epsilon}$ (resp.
 $\overleftarrow{P_\epsilon}$)
of $P$ is defined by:
$\overrightarrow{P_\epsilon}=\{Q\in \F \ \text{s.t.} \ \ \epsilon_{PQ}=1\}$  (resp.
$\overleftarrow{P_\epsilon}=\{Q\in \F \ \text{s.t.} \ \ \epsilon_{PQ}=-1\}$). 
\item  Let ${\mathbb O_\F}$ be the set of $\FF-$valued functions  on the Fano cube $V_\F$.
\end{enumerate}
\end{definition}

\noi
Recall that a trichotomous binary relation on $\F$ is a binary relation  $R$ such that for all $x,y \in \F$ exactly one of the following holds: $x=y$ or $x R y$ or $y R x$. A multiplication factor $\epsilon$
uniquely determines
such a relation on $\F$ by: $x R y$ {\it iff} $\epsilon_{xy}=1$. Conversely, every trichotomous binary relation on $\F$ is obtained from a unique multiplication factor in this way so that trichotomous binary relations and multiplication factors are the same thing. 
\\

In the rest of the paper an important r\^ole will played by the ``exponential'' map $e: \mathbb Z_2 \to \{-1,1\}$ and its inverse $\ell : \{-1,1\} \to \mathbb Z_2$  defined  by
\beqa
 e^0=1\ \  , \ \ e^1=-1 \ \ 
\text{and} \ \
\ell(1)=0 \ \ , \ \ \ell(-1)=1 \ .\nn
\nn
\eeqa

\begin{remark}
It is easy to see that if $N$ is  a norm  then the set $\big\{P \in \F \ \mbox{s.t.}\  N(P)=1\big\}$ is either $\F$ or a line. This means there is a unique linear form $n \in V_\F^\ast$ such
that $N(P)= e^{n(P)}$ for all $P \in\F$.
\end{remark}

For $P\in V_\F$  define $e_P\in {\mathbb O_\F}$ by
\beqa
e_P(Q) = \left\{
\begin{array}{ll}
1& \mathrm{~if~} Q=P\\
0& \mathrm{~if~} Q\ne P\ . 
\end{array}\right.\nn
\eeqa
Then $\{e_P \ \text{s.t.} \ P\in V_\F\}$ is an $\FF-$basis of ${\mathbb O_\F}$ and
\beqa
{\mathbb O_\F} = \FF e_0 \oplus \mathrm{Vect}\Big<e_P\ \text{s.t.}\  P \in \F\Big> \ . \nn
\eeqa
A norm  and a multiplication factor on $\F$ allow us to endow ${\mathbb O_\F} $ with a norm and a  multiplication:

\begin{definition}\label{def:mo}
Let $\F$ be a Fano plane equipped with a norm $N$ and a multiplication factor $\epsilon$.
\begin{enumerate}
\item
The multiplication $\cdot_\e : {\mathbb O_\F} \times {\mathbb O_\F} \to {\mathbb O_\F}$ is the unique bilinear map such that
\begin{enumerate}
\item For all $P\ne Q \in \F$:\ \ \ \  $e_P\cdot_\e e_Q = \epsilon_{PQ} e_{P + Q}$;
\item For all $P \in \F$:\ \ \ \   $e_P\cdot_\e e_P = -N(P)e_0$;
\item For all $P \in V_\F  $:\ \ \ \   $e_0 \cdot_\e e_P = e_P \cdot_\e e_0 = e_P$.
\end{enumerate}
\item The norm $N_{\mathbb O_\F}: {\mathbb O_\F} \to \mathbb F$ is the  quadratic form :
$N_{\mathbb O_\F} (\lambda^0 e_0 + \sum \limits_{P\in\F} \lambda^P e_{P})= (\lambda^0)^2 + \sum\limits_{P\in \F} (\lambda^P)^2 N(P)$.

\end{enumerate}
\noindent We denote by ${\bf 1}$ the quadratic form on ${\mathbb O_\F}$ associated to the trivial norm on $\F$. By definition
 the triple  $({\mathbb O_\F},N_{\mathbb O_\F},\e )$ is a composition algebra  iff $N_{\mathbb O_\F} (Z  \cdot_\e W) =N_{\mathbb O_\F}(Z) N_{\mathbb O_\F}(W), \forall Z, W \in {\mathbb O_\F}$.
\end{definition}

\begin{example}
A multiplication factor  on $\F$ can be deduced from the system of arrows  in Figure \ref{fig:oct} left panel as follows: given distinct $P, Q$ in $ \F$,
 set $\epsilon_{PQ}=1$ {\it iff} there is an arrow from $P$ to $Q$. Then one can show that 
 $({\mathbb O_\F}, {\bf 1}, \epsilon)$ is a composition algebra. This is Freudenthal's original construction.
\end{example}

In view of (c) above,  we denote $e_0$ by $1$  and  now introduce the ``compositor'' which measures the failure of  $({\mathbb O_\F},N_{\mathbb O_\F},\e )$ to be a composition algebra.
\begin{definition}
The compositor $C: {\mathbb O_\F} \times {\mathbb O_\F} \to \mathbb F$ is defined by
\beqa C(Z,W)= N_{\mathbb O_\F}(Z \cdot_\e W) - N_{\mathbb O_\F}(Z) N_{\mathbb O_\F}(W) \ \ \forall Z, W \in {\mathbb O_\F} \
 .\nn \eeqa
\end{definition}

\noi
The compositor can be expressed as a sum over all lines and quadrilaterals:

\begin{proposition}\label{prop:comp}
With the notation above, let $Z=\lambda^0 + \sum \limits_{P \in \F} \lambda^P e_P \in {\mathbb O_\F}$  and let $W=\mu^0 + \sum \limits_{P\in \F} \mu^P e_P \in {\mathbb O_\F}$, where $\lambda^0,\mu^0, \lambda^P, \mu^P \in \mathbb F$. Then the compositor
is given by
\beqa
C(Z,W)= &\hskip -.3truecm \sum\limits_{(P,Q,R,S) \in \mathrm{Quadrilateral}} N(P+ Q) \epsilon_{PQ} \epsilon_{RS} \lambda^P \mu^Q \lambda^R \mu^S\nn\\
&\hskip -1.2 truecm +
\sum\limits_{(P,Q,R) \in \mathrm{Line}} 2 N(R)\epsilon_{PQ} \lambda^P \mu^Q\big(\lambda^0 \mu^R + \lambda^R \mu^0 \big) \ . \nn
\eeqa
 In particular $({\mathbb O_\F}, N_{\mathbb O_\F}, \epsilon)$ is a composition algebra  {\it iff}:
\beqa
\label{eq:comp1}
(i)&&N(P+R) \epsilon_{PQ} \epsilon_{QR} = 1\hskip 4.1truecm \text{for any  line }\{P,Q,R\},\\
\label{eq:comp2}
(ii)&& N(P+Q)\;\epsilon_{PQ}\epsilon_{QR}\epsilon_{RS}\epsilon_{SP} \;N(P+S) =-1\quad \text{for any  quadrilateral }\{P,Q,R,S\}.
\eeqa
\end{proposition}
\begin{demo}
Let  $\F = \{P_1,\cdots, P_7\}$. 
By calculation we obtain
\beqa
N_{\mathbb O_\F}(Z) N_{\mathbb O_\F}(W)= \Big((\lambda^0)^2 + \sum_{i=1}^7 (\lambda^{P_i})^2 N(P_i) \Big)
\Big( (\mu^0)^2 + \sum_{i=1}^7 (\mu^{P_i})^2 N(Q_i) \Big) \  ,\nn
\eeqa
and
\beqa
N_{\mathbb O_\F}(Z \cdot_\e W) &=&\Bigg(\lambda^0 \mu^0 - \sum\limits_{i=1}^7 \lambda^{P_i} \mu^{P_i} N(P_i) \Bigg)^2 \nn\\
&&+ \sum \limits_{i=1}^7
\Bigg(\lambda^0 \mu^{P_i} + \lambda^{P_i} \mu^0 + \sum\limits_{\tiny \begin{array}{c} R,S \\ \text{s.t.} \\ R+S= P_i \end{array}}\epsilon_{RS} \lambda^R \mu^S\Bigg)^2 N(P_i) \ . \nn
\eeqa
Taking the difference  of these two expressions, after some simplifications we obtain the desired expression for the compositor $C(Z,W)$.
\end{demo}

If $({\mathbb O_\F}, {\bf 1},\epsilon)$ is a composition algebra, we have the following `Fano geometric' consequences.
 Later on we will show   converses to this theorem (see Theorems \ref{theo:future1} and \ref{theo:future2}).
\begin{theorem} \label{theo:future}
Let $({\mathbb O_\F}, {\bf 1}, \epsilon)$ be a composition algebra. Then 
  either  $\overrightarrow{P}$   is a line for all $P$ in $\F$, 
or  $\overleftarrow{P}$  is a line for all $P$ in $\F$. In the first case the set of triangles $\{\overleftarrow{P} \ \ \mbox{s.t.} \ \ P\in\F\}$ satisfies the axioms of a Fano plane  and in the second the set of triangles $\{\overrightarrow{P} \ \ \mbox{s.t.} \ \ P\in\F\}$ satisfies the axioms of a Fano plane.

\end{theorem}
\begin{demo}
We begin by proving a series of lemmas.
\begin{lemma}
Let ${\cal Q}$ be a quadrilateral in $\F$. Then  there exists  a unique triangle ${\cal T}=\{P,Q,R\}$ contained in ${\cal Q}$ which is oriented, {\it i.e.}, such that $\epsilon_{PQ}=\epsilon_{QR}=\epsilon_{RP}$. If ${\cal Q}={\cal T}\cup\{S\}$ then either $\overrightarrow{S}={\cal T}$ or $\overleftarrow{S}={\cal T}$.
\end{lemma}
\begin{demo}
Suppose for contradiction   that the four triangles in  ${\cal Q}$ are not oriented. Let ${\cal Q}= \{P,Q,R,S\}$ and consider the triangle $\{P,Q,R\}$. Without loss of generality we can assume
$$\epsilon_{PQ}=1, \ \ \epsilon_{QR}=1 \ \ \mbox{and} \ \ \epsilon_{RP}=-1 \ . $$

If $\epsilon_{QS}=1$ then  from the non-orientability of the triangle $\{P,Q,S\}$ we must have $\epsilon_{SP}=-1$. Hence
\beqa
\epsilon_{PR} \epsilon_{RQ} \epsilon_{QS} \epsilon_{SP}=(1)(-1)(1)(-1)=1 \nn
\eeqa
which is in contradiction with \eqref{eq:comp2}. If $\epsilon_{QS}=-1$ then  from the non-orientability of the triangle $\{Q,R,S\}$ we must have $\epsilon_{RS}=-1$ and
\beqa
\epsilon_{PR} \epsilon_{RS} \epsilon_{SQ} \epsilon_{QP}=(1)(-1)(1)(-1)=1 \nn
\eeqa
which is also in contradiction with \eqref{eq:comp2}. Hence  there exists a least one oriented triangle, say $\{P,Q,R\}$, in ${\cal Q}$. By \eqref{eq:comp2}
\beqa
\epsilon_{SP} \epsilon_{PR} \epsilon_{RQ} \epsilon_{QS} = -1 \nn
\eeqa
and this implies that
\beqa
\label{eq:sp}
\epsilon_{SP}  \epsilon_{SQ} = 1 
\eeqa
since $\{P,Q,R\}$ is oriented. Hence  $\epsilon_{SP} = \epsilon_{SQ}$ and similarly $\epsilon_{SP} = \epsilon_{SR}$.
This means that no triangle in ${\cal Q}$ containing $S$ is oriented and, by \eqref{eq:comp1}, that either $\overrightarrow{S}={\cal T}$ or $\overleftarrow{S}={\cal T}$.
\end{demo}

The above lemma means that every quadrilateral has a distinguished point:  the complement of its oriented triangle which is also the orthopoint of its oriented triangle ({\it c.f.} L7 in List \ref{li:inci}).
\begin{lemma}
Let ${\cal Q}_1$ and ${\cal Q}_2$ be quadrilaterals in $\F$ and let $P_1$ and $P_2$ be the associated distinguished points. Then  ${\cal Q}_1 ={\cal Q}_2$ {\it iff} $P_1 =P_2$. 
\end{lemma}
\begin{demo}
Suppose for contradiction that ${\cal Q}_1$ and ${\cal Q}_2$ are distinct quadrilaterals with the same distinguished point $P$.
Since ${\cal Q}_1 \cap {\cal Q}_2$ contains two points ({\it c.f.} L4 in List \ref{li:inci})  we can choose $R_1$ in $ {\cal Q}_1 \setminus {\cal Q}_2$. Since $R_1\ne P$
there is a unique point $R_2\in {\cal Q}_2$ such that $\{R_1,P,R_2\}$ is a line. Hence by  \eqref{eq:comp1} $\epsilon_{PR_1} \epsilon_{PR_2}=-1$ which contradicts \eqref{eq:sp}.
This proves the lemma.
\end{demo}
\begin{lemma}\label{lem:trintsect}
Let  ${\cal Q}_1$ and ${\cal Q}_2$ be distinct quadrilaterals, let ${\cal T}_1, {\cal T}_2$ be the corresponding oriented triangles and
let $P_1$, $P_2$  be the corresponding distinguished  points.
The following are equivalent
\begin{enumerate}
\item[(i)] ${\cal T}_1$ and ${\cal T}_2$ intersect in one point.

\item[(ii)] ${\cal Q}_1 \cap {\cal Q}_2$   contains  exactly one of $P_1$ or $P_2$.
\end{enumerate}
\end{lemma}
\begin{demo}
(i) $\Rightarrow$ (ii): 
Let  ${\cal T}_1 \cap {\cal T}_2=\{R\}$ and  ${\cal Q}_1 \cap {\cal Q}_2=\{R,S\}$.
If  $S \not \in {\cal T}_1 \cup {\cal T}_2$, 
then ${\cal Q}_1 = {\cal T}_1 \cup \{S\}$ and  ${\cal Q}_2 = {\cal T}_2 \cup \{S\}$ and hence $S$ is the distinguished
point of ${\cal Q}_1$ and ${\cal Q}_2$ which contradicts the previous lemma.  This shows that one of the triangles, say ${\cal T}_1$ contains $S$ and that $S$ is the distinguished point of 
 ${\cal T}_2$.

(i) $\Leftarrow$ (ii): Since by hypothesis either $P_1$ or $P_2$ is in  ${\cal Q}_1 \cap {\cal Q}_2$, without loss of generality we can suppose that ${\cal Q}_1 \cap {\cal Q}_2=\{P_1,R\}$ where
$R \in {\cal T}_2$. However $R\in {\cal T}_1$ since $R \neq P_1$ and so  $R \in {\cal T}_1 \cap {\cal T}_2$. Moreover  ${\cal T}_1 \cap {\cal T}_2 \subset {\cal Q}_1 \cap {\cal Q}_2 = \{P_1,R\}$ and
$P_1 \not \in {\cal T}_1$. Hence ${\cal T}_1 \cap {\cal T}_2= \{R\}$.
\end{demo}
\begin{lemma}
Let ${\cal Q}_1$ and ${\cal Q}_2$ be two quadrilaterals with distinguished points $P_1,P_2$ respectively.   Then either $P_1$ or $P_2$ is in  ${\cal Q}_1 \cap {\cal Q}_2$.
\end{lemma}
\begin{demo}
Let ${\cal Q}_1 \cap {\cal Q}_2 = \{R,S\}$ and suppose for contradiction that $\{P_1,P_2\}\cap\{ R,S\}=\emptyset$. Define $Q_1 \in {\cal Q}_1$ and $Q_2 \in {\cal Q}_2$ by
${\cal Q}_1=\{P_1,Q_1,R,S\}$ and ${\cal Q}_2=\{P_2,Q_2,R,S\}$.
Since $\{R,Q_1,S\}$ is a triangle, the points of $\F$ are: $R,Q_1, S, R+S, S+Q_1, Q_1+R$ and $R+S+Q_1=P_1$. It follows that  ${\cal Q}_2=\{R,S,Q_1+R,Q_1+S\}$ as $R+Q_1, Q_1+S, S+R$ is the ortholine of
  $\{R,Q_1,S\}$. Hence either $P_2=Q_1+R$ or $P_2=Q_1+S$. In the first case we have $Q_2=Q_1+S$ and
  \beqa
P_1+P_2= S=Q_1+Q_2 \nn
\eeqa
which means that the line through $P_1, P_2$ intersects the line through $Q_1,Q_2$ in $S$.

By \eqref{eq:comp2} for the quadrilaterals ${\cal Q}_i$
\beqa
\epsilon_{RS} \epsilon_{SQ_i} \epsilon_{Q_iP_i} \epsilon_{P_iR}=-1, \ \ i=1,2 \nn 
\eeqa
which  using \eqref{eq:sp} simplifies to
\beqa
 \epsilon_{SQ_1} =  \epsilon_{SQ_2}=1 \nn
\eeqa
and this contradicts \eqref{eq:comp1} for the line $\{Q_1,S,Q_2\}$.
\end{demo}

\begin{lemma}\label{lem: trinsect1}
Let ${\cal Q}_1$ and ${\cal Q}_2$ be two quadrilaterals with distinguished points $P_1,P_2$ respectively.   Then exactly one of  $P_1$ or $P_2$ is in  ${\cal Q}_1 \cap {\cal Q}_2$.
\end{lemma}
\begin{demo}
Suppose for contradiction that  ${\cal Q}_1 \cap {\cal Q}_2 = \{P_1,P_2\}$. Define $R$ and $S$ in ${\cal Q}_1$ by ${\cal Q}_1=\{R,S,P_1,P_2\}$.  Since $\{R,S,P_1\}$ is a triangle
${\cal Q}_2=\{P_1,P_2=R+S+P_1,P_1+S,P_1+R\}$. By \eqref{eq:comp1} for the lines $\{S+P_1,P_2,R\}$ and   $\{R+ P_1,P_2,S\}$ we have
\beqa
\epsilon_{P_2R}=-\epsilon_{P_2S+P_1} \ \ \mbox{and} \ \ \epsilon_{P_2S}=-\epsilon_{P_2R+P_1} \nn
\eeqa
However from \eqref{eq:sp} since $P_2$ is the distinguished point of ${\cal Q}_2$ it follows that $\epsilon_{P_2R}=\epsilon_{P_2S}$ which contradicts the fact that the triangle
$\{P_2,S,R\}$ is oriented.
This proves the lemma.
\end{demo}
\begin{lemma}
Let  ${\cal Q}_1$ and ${\cal Q}_2$ be distinct quadrilaterals, let ${\cal T}_1, {\cal T}_2$ be the corresponding oriented triangles and $P_1$, $P_2$ the corresponding distinguished points.
Then
\beqa
\epsilon_{P_1 R_1} = \epsilon_{P_2 R_2} \ \ \forall (R_1,R_2) \in {\cal T}_1\times {\cal T}_2 \nn \ .
\eeqa
\end{lemma}
\begin{demo}
Without loss of generality by  the lemmas above  we can suppose that  $P_1 \in {\cal Q}_1 \cap {\cal Q}_2$ and $P_2 \not \in {\cal Q}_1 \cap {\cal Q}_2$.
The third point  $Q$ on the  line  through  $P_1 P_2$ is in ${\cal Q}_1$. Thus by \eqref{eq:comp1} $\epsilon_{P_1Q} = \epsilon_{P_2 P_1}$ and the lemma follows from \eqref{eq:sp}.
\end{demo}

\noi
{\it Proof of the Theorem:}
Let ${\cal Q}_i, i=1,\cdots,7$ be the quadrilaterals of $\F$,  let ${\cal T}_i,i=1,\cdots,7$ be the corresponding oriented triangles and let $P_i, i=1,\cdots, 7$ be the corresponding distinguished points.
 By the above lemmas the composition factor $\epsilon$ has  one of the following properties
 \beqa
\epsilon_{P_i R} =1\ \ \forall (P_i,R)\in\F\times{\cal T}_i \quad\mbox{or} \quad
\epsilon_{P_i R} =-1 \ \  \forall (P_i,R)\in\F\times{\cal T}_i \ . \nn
\eeqa
In the first case we say that $\epsilon$ is a ``source'' and in the second a ``sink''.

If $\epsilon$ is a source (resp. sink) then for all $P_i$  we have
$\overrightarrow{P_i}={\cal T}_i$ (resp. $\overleftarrow{P_i}={\cal T}_i$)
is a triangle whose orthopoint is $P_i$ and hence $\overleftarrow{P_i} $ (resp. $\overrightarrow{P_i})$
is a line ({\it c.f.} L2 in  List \eqref{li:inci}).  This proves the first part of of the theorem.

We only prove the second part of the theorem in the case where $\overrightarrow{P_i}$ is a line for all $P_i$.
Since $\overleftarrow{P}_i={\cal T}_i$ (see above), it follows from    Lemmas \ref{lem:trintsect} and \ref{lem: trinsect1} that the set of triangles $\{\overleftarrow{P_i} \ \ \mbox{s.t.} \ \ P_i\in\F\}$ satisfies axiom  \ref{def: Fano plane} (b).  To prove axiom \ref{def: Fano plane} (a) we have to show that given  distinct points $P_i,P_j$ there exists a unique $\overleftarrow{P_k} $ containing them. This follows from the fact that  $P_i,P_j\in\overleftarrow{P_k} $ iff $P_k\in\overrightarrow{P_i}\cap\overrightarrow{P_j}$, and the fact that the two lines $\overrightarrow{P_i},\overrightarrow{P_j}$ intersect in a unique point. 
\end{demo}
\begin{definition}\label{def:Opm}
We denote ${\mathbb O_\F}_1^+$  (resp.  ${\mathbb O_\F}_1^-$)
the set of all composition algebras $({\mathbb O_\F}, {\bf 1},\epsilon)$ such that 
$\overrightarrow{P}$ (resp. $\overleftarrow{P}$)  is a line for all $P\in\F$.
\end{definition}

In fact the condition to be a composition algebra is  implied by either of the other two conditions in this definition
(see Theorem \ref{theo:future1}). 
\section{Oriented maps and $n-$oriented maps}   

Let $n$ be a linear form on $V_\F$.
In this section we introduce and classify $n-$oriented maps from $\F \to \F^\ast$ (if $n \equiv 0$ we refer to these maps simply as oriented maps). These can be thought of as  logarithms of the composition factors of the previous section with the case $n\equiv0$ corresponding to division  composition algebras.

\begin{definition}\label{def:contra}
An oriented map is a map $\alpha : \F \to \F^\ast$ such that
\begin{enumerate}
\item[(i)] For all $P\in \F, \;P \not \in \alpha(P)\ .$
\item[(ii)] If $P\ne Q \in \F$ then: \; $P \in \alpha(Q) \Leftrightarrow Q \not\in \alpha(P)\ .$
\end{enumerate}
\end{definition}

Since $P \not\in \alpha(P),$ the triple $\Delta_P=\F\setminus\{P\} \setminus \{\alpha(P)\}$ is a triangle and $P$ is the orthopoint
of $\Delta_P$ ({\it cf.}\ L3 in List \ref{li:inci}). This shows that the map $\alpha: \F \to \F^\ast$ is a bijection.

\begin{example}
An oriented map can be deduced from the system of arrows in Figure \ref{fig:oct} left panel: given $P$ in $\F$ define
$\alpha(P)= \overrightarrow{P}$. 
\end{example}

\begin{proposition} \label{prop:trian}
\begin{enumerate}
\item[(i)]
Let $\alpha: \F \to \F^\ast$ be an oriented map. For all $P\in \F$ let $\Delta_P= \F \setminus\{P\} \setminus \alpha(P)$.
Then the set of seven triangles $\Delta_\alpha=\{\Delta_P \ \ \text{s.t.} \ \  P \in\F\}$ satisfies the axioms of a Fano plane:
(a)  two distinct points of $\F$  are contained in a unique element of $\Delta_\alpha$ and (b) two distinct elements of $\Delta_\alpha$ intersect in a unique point.
\item[(ii)] Conversely, if $\Delta$ is a set of seven triangles satisfying the axioms
of the Fano plane there  exists a unique oriented map $\alpha:\F\to \F^\ast$ such that
$\Delta=\Delta_\alpha$.
\end{enumerate}
\end{proposition}

\begin{demo}
(i)(a):  This is immediate since if $P\ne Q \in\F$ and $R \in \F$, then by  \ref{def:contra} (ii):
\beqa
  R \in \alpha(P) \cap \alpha(Q)\ \  \Leftrightarrow \ \ P, Q \in\Delta_R \ . \nn
\eeqa

(i)(b):  By definition, for all $P,Q$ in $\F$ we have $\Delta_P\cap \Delta_Q= \Big(P\cup \alpha(P)\cup Q\cup\alpha(Q)\Big)^c$. If $P\ne Q$ then we can suppose without loss
of generality that $P\in\alpha(Q)$ and $Q\in\Delta_P$ ({\it cf.} L2 in List \ref{li:inci}). Hence 
$\Big(P\cup \alpha(P)\Big)\cap \Big(Q\cup \alpha(Q)\Big)=P\cup\Big(\alpha(P)\cap \alpha(Q)\Big)$ and
\beqa
\mbox{Card}\Big[\Big(P\cup \alpha(P)\Big)\cup \Big(Q\cup \alpha(Q)\Big)\Big]=6\ .\nn
\eeqa
This proves the result.

(ii): Let $\Delta=\{\Delta_1,\cdots, \Delta_7\}$. Then we have 7 partitions
\beqa
\F = \Delta_i \cup L_i \cup P_i  \ \ i=1,\cdots,7\nn 
\eeqa
where $L_i$ is the ortholine of $\Delta_i$ and $P_i$ is the orthopoint.

If $L_i=L_j$ then $\Delta_i,\Delta_j$ are two triangles in the same quadrilateral $\F \setminus L_i$ which intersect in either three points
or two points. Since $\Delta$ satisfies the axioms of the Fano plane we must have $\Delta_i=\Delta_j$. This shows that the map
$\Lambda:\Delta_i
\mapsto L_i$ is a bijection of $\Delta$ with $\F^\ast$.

If $i\ne j$ suppose for contradiction that  $P_i=P_j$. Then, since the triangles $\Delta_i, \Delta_j$ are distinct and  the distinct quadrilaterals $Q_i=\F\setminus L_i$ and
$Q_j=\F\setminus L_j$ intersect in two points ({\it cf.} L4 in List \ref{li:inci}) we have

\beqa
Q_i \cap Q_j = \big\{P_i\big\} \cup \big(\Delta_i \cap \Delta_j\big) \ .\nn
\eeqa
Let $\Delta_k$ be the unique triangle such that its ortholine $L_k$ contains $P_i$ and $\Delta_i \cap \Delta_j$. Then the quadrilateral
$Q_k=\F\setminus L_k$ is distinct from $Q_i$ and $Q_j$ and satisfies $Q_i \cap Q_j \cap Q_k = \varnothing$. Hence ({\it cf.}
L4 in \ref{li:inci})
\beqa
Q_k&=& \big(Q_k \cap Q_i\big) \cup\big(Q_k \cap Q_j\big) =\big(Q_k \cap \Delta_i\big) \cup\big(Q_k \cap \Delta_j\big)\nn\\
&=&\big(\{P_k\}\cap \Delta_i\big) \cup \big(\Delta_k \cap \Delta_i\big)\cup \big(\{P_k\}\cap \Delta_j\big)\cup  \big (\Delta_k \cap \Delta_j\big)\
\  \nn
\eeqa
and therefore either $\Delta_k \cap \Delta_i$ or  $\Delta_k \cap \Delta_j$ contains two points. This is a contradiction and
the map $\Pi:\Delta_i \mapsto  P_i$ is a bijection of  $\Delta$ with  $\F$.\\

Define $\alpha: \F \to \F^\ast$ by $\alpha = \Lambda \circ \Pi^{-1}$.
 We now show that $\alpha$ is an oriented map. By construction for all $P\in \F$ we have
\beqa
\F= \{P\} \cup \alpha(P) \cup \Lambda^{-1} \circ \alpha(P) \nn
\eeqa
and hence $P \not \in \alpha(P)$. This proves (i) of Definition \ref{def:contra}.

Let $P\ne Q \in \F$.  Then
\beqa
\F &=&  \{P\} \cup \alpha(P) \cup \Lambda^{-1} \circ \alpha(P) \nn\\
 &=&  \{Q\} \cup \alpha(Q) \cup \Lambda^{-1} \circ \alpha(Q) \nn\ .
\eeqa
Suppose for contradiction that $P \in \alpha(Q)$ and  $Q \in \alpha(P)$. The union $L_{PQ}$
 of the  distinct lines $\alpha(P), \alpha(Q)$
 contains five points and its complement $L_{PQ}^c$  two.
Since $\{P\} \cup \alpha(P)\subset L_{PQ}$, $L^c_{PQ}\subset \Lambda^{-1}\circ \alpha(P)$.
Similarly   $L^c_{PQ}\subset \Lambda^{-1}\circ \alpha(Q)$ and hence $\Lambda^{-1}\circ \alpha(P) \cap 
\Lambda^{-1}\circ \alpha(Q)$ contains at least two points
which, since $\Lambda^{-1} \circ \alpha (P) \ne\Lambda^{-1} \circ \alpha (Q)$, is a contradiction to the second Fano axiom (see Definition \ref{def:contra}). This  proves (ii) of Definition \ref{def:contra}.
\end{demo}

In order to account for split octonions we introduce the following definition which generalises Definition \ref{def:contra}:

\begin{definition}\label{def:contra-gene}
Let $n \in V_\F^\ast$.  An $n-$oriented map is a map $\alpha: \F \to V_\F^\ast$ such that (writing $\alpha_P$ for $\alpha(P)$):
\begin{itemize}
\item[(i)] For all $P\in \F, \ \ \alpha_P(P) \neq 0 \Leftrightarrow P \in \mbox{Ker}(n)$ \ \ \ (equivalently $\alpha_P(P) + n(P)=1$)\ .
\item[(ii)]  If $P\not = Q \in \F$ then \ \ $\alpha_P(Q)+ \alpha_Q(P) =1$ \ . 
 \end{itemize}
 The triple $(\F,n, \alpha)$ is called an $n-$oriented  Fano plane and we denote by $\F_n$ the set of all $n-$oriented Fano planes.
 Note that $n \not\in $  \rm{Im}$(\alpha)$ by (i).
\end{definition}

Identifying  a line in $\F$ with the unique  linear form on $V_\F$ of which it is the kernel, 
it is clear that a $0-$oriented map  is the same thing as an oriented  map in the sense of Definition \ref{def:contra}. In fact as we now show there is a natural bijection between $0-$oriented maps
and $n-$oriented maps for $n \in V_\F^\ast$.

\begin{proposition}\label{prop:corres}
Let $\alpha$ be an  oriented  map and let $n \in V_\F^\ast$. Define $\alpha': \F \to V_\F^\ast$ by
\beqa
\alpha'_P(Q) = \alpha_P(Q) + n(P) n(Q) \ .\nn
\eeqa
Then $\alpha'$ is an $n-$oriented map and every $n-$oriented  map is uniquely obtained in this way.
In particular $\alpha \mapsto \alpha'$  is a natural bijection $\F_0 \cong \F_n$.
\end{proposition}

\begin{demo}
We have to show that
\begin{itemize}
\item[(i)] $\alpha'_P(P) \neq 0 \Leftrightarrow P \in \mbox{Ker}(n)$;
\item[(ii)] $\alpha'_P(Q)+ \alpha'_Q(P) =1 \ \ P\neq Q$.
 \end{itemize}
 The second property is evident. For the first property:  if $n(P)=0$ then $\alpha'_P(P) =\alpha_P(P) =1$ whereas if $n(P)=1$ then $\alpha'_P(P) =\alpha_P(P) +n(P)=0$.
 Conversely, if $\alpha'$ is an $n-$oriented map the formul\ae\ above define a $0-$oriented map $\alpha$.
\end{demo}

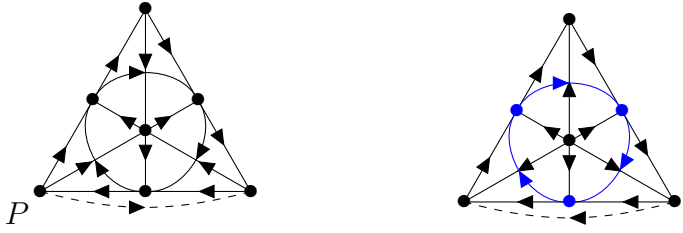
\begin{figure}[!ht]
  \begin{center}\hskip 2.truecm
  \begin{minipage}{4cm}
    \begin{tikzpicture}[scale=.7]
\tikzstyle{point}=[circle,draw]

\tikzstyle{ligne}=[thick]
\tikzstyle{pointille}=[thick,dotted]

\tikzset{->-/.style={decoration={
  markings,
  mark=at position .5 with {\arrow{>}}},postaction={decorate}}}

\tikzset{middlearrow/.style={
        decoration={markings,
            mark= at position 0.5 with {\arrow{#1}} ,
        },
        postaction={decorate}
    }}

\coordinate  (3) at ( -2,-1.15);

\coordinate  (2) at ( 0,-1.15) ;
\coordinate  (5) at ( 2,-1.15);

\coordinate  (6) at (0,2.31);
\coordinate  (1) at (1, .58);

\coordinate  (4) at (-1., .58) ;
\coordinate  (7) at (0,0) ;

\draw [color=black, middlearrow={triangle 45}] (1) to [bend left=65] (2);
\draw [color=black,  middlearrow={triangle 45}] (2) --(3);
\draw [color=black, middlearrow={triangle 45}] (3) --(4);
\draw [color=black,  middlearrow={triangle 45}] (6) --(7);
\draw [color=black,  middlearrow={triangle 45}] (7) --(1);
\draw[dashed,middlearrow={triangle 45}] (3) to [bend right=15]  (5);

\draw [middlearrow={triangle 45}] (6) --(1);
\draw [ middlearrow={ triangle 45}] (1) --(5);

\draw[middlearrow={triangle 45}]  (5) --(2);

\draw[middlearrow={triangle 45}]  (3) --(7);

\draw  [color =black, middlearrow={triangle 45}]  (5) --(7);
\draw [color =black, middlearrow={triangle 45}] (7) --(4);
\draw[middlearrow={triangle 45}]  (4) --(6);
\draw [middlearrow={triangle 45}]  (7) --(2);

\draw [middlearrow={triangle 45}] (2) to [bend left=65] (4);
\draw [middlearrow={triangle 45}]  (4) to [bend left=65] (1);
\draw (1)  node[scale=0.6] {$\bullet$} ;
\draw (1)  node {$\bullet$} ;

\draw (2)  node{$\bullet$} ;
\draw (3) node[below left=0pt]{$P$} ;
\draw (3) [color=black]  node{$\bullet$} ;
\draw (4)  node{$\bullet$} ;
\draw (5)  node{$\bullet$} ;
\draw (6)  node{$\bullet$} ;
\draw (7)  node{$\bullet$} ;

  \end{tikzpicture} 
  \end{minipage}\hskip 2.truecm
  \begin{minipage}{8cm}
 \begin{tikzpicture}[scale=.7]
\tikzstyle{point}=[circle,draw]

\tikzstyle{ligne}=[thick]
\tikzstyle{pointille}=[thick,dotted]

\tikzset{->-/.style={decoration={
  markings,
  mark=at position .5 with {\arrow{>}}},postaction={decorate}}}

\tikzset{middlearrow/.style={
        decoration={markings,
            mark= at position 0.5 with {\arrow{#1}} ,
        },
        postaction={decorate}
    }}

\coordinate  (3) at ( -2,-1.15);

\coordinate  (2) at ( 0,-1.15) ;
\coordinate  (5) at ( 2,-1.15);

\coordinate  (6) at (0,2.31);
\coordinate  (1) at (1, .58);

\coordinate  (4) at (-1., .58) ;
\coordinate  (7) at (0,0) ;

\draw [color=blue, middlearrow={triangle 45}] (1) to [bend left=65] (2);
\draw [color=black,  middlearrow={triangle 45}] (2) --(3);
\draw [color=black, middlearrow={triangle 45}] (3) --(4);
\draw [color=black,  middlearrow={triangle 45}] (7) --(6);
\draw [color=black,  middlearrow={triangle 45}] (7) --(1);
\draw[dashed,middlearrow={triangle 45}] (5) to [bend left=15]  (3);

\draw [middlearrow={triangle 45}] (6) --(1);
\draw [ middlearrow={ triangle 45}] (1) --(5);

\draw[middlearrow={triangle 45}]  (5) --(2);

\draw[middlearrow={triangle 45}]  (7) --(3);

\draw  [color =black, middlearrow={triangle 45}]  (7) --(5);
\draw [color =black, middlearrow={triangle 45}] (7) --(4);
\draw[middlearrow={triangle 45}]  (4) --(6);
\draw [middlearrow={triangle 45}]  (7) --(2);

\draw [color=blue, middlearrow={triangle 45}] (2) to [bend left=65] (4);
\draw [color=blue,middlearrow={triangle 45}]  (4) to [bend left=65] (1);
\draw (1)  node[scale=0.6] {$\bullet$} ;
\draw (1)  node[color=blue] {$\bullet$} ;

\draw (2)  node[color=blue]{$\bullet$} ;
\draw (3) [color=black]  node{$\bullet$} ;
\draw (4)  node[color=blue]{$\bullet$} ;
\draw (5)  node{$\bullet$} ;
\draw (6)  node{$\bullet$} ;
\draw (7)  node{$\bullet$} ;

  \end{tikzpicture}
\end{minipage}
  
  \caption{
From oriented maps to $n-$oriented maps: change direction of arrows on the six line segments containing no blue point (see Proposition \ref{prop:corres}).}
\label{fig:contr}
\end{center}
\end{figure}

\begin{example}
Let $n$ be the linear form whose  kernel is the blue line in Figure \ref{fig:contr}, right panel. An $n-$oriented map $\alpha$ can be deduced from the system of arrows in the right Fano plane 
as follows: for every point $Q$ except the central point,  define $\alpha'(Q)$ to be the unique line containing
$\overrightarrow{Q}$  and for $Q_0$ the central point, define
 $\alpha'(Q_0)$ to be $0$. This is the $n-$oriented  map  $\alpha'$ obtained from the $0-$oriented map $\alpha$ corresponding to  Figure \ref{fig:contr}, left panel by Proposition \ref{prop:corres}.
Observe that $Q_0$ is in the ``past'' of every other point.
\end{example}

 Let  ${\cal S}_0^2(V_\F^\ast)$ be  the space of  bilinear  forms $B$ on $V_\F$ satisfying $B(P,P)=0$ for all $P$ in $\F$.
We now show that  $\F_0$, the space of oriented maps from $\F$ to $V_{\F}^\ast$, is an affine space for ${\cal S}_0^2(V_\F^\ast)$.

\begin{proposition}
\begin{enumerate}
\item[(i)] Let $P\mapsto \alpha_P$ and $P\mapsto  \beta_P$ be in $\F_0$, and define $\Delta: V_\F \times V_\F \to \mathbb Z_2$ by
\beqa
\Delta(P,Q) =\left\{\begin{array}{cl} \alpha_P(Q) + \beta_P(Q) & \forall P,Q \in \F \times \F\ ,\\
                                      0&\mbox{if} \ P  \ \mbox{or} \ Q=0\ .
                                      \end{array}\right.\nn
                                      \eeqa
Then $\Delta$ is symmetric, bilinear and $\Delta(P,P)=0$ for all $P$ in $V_\F$.
\item[(ii)] Let $P \mapsto \alpha_P$ be in $\F_0$ and $\Delta:  V_\F \times V_\F \to \mathbb Z_2$ be symmetric, bilinear and such that $\Delta(P,P)=0$ for all $P$ in $\F$. Define
$\beta_P$ by 
\beqa
\beta_P(Q)= \Delta(P,Q) + \alpha_P(Q)  \ , \ \ \forall Q \in V_\F\ . \nn
\eeqa
Then $P \mapsto \beta_P$ is an element of $\F_0$.
\end{enumerate}  
\end{proposition}
\begin{demo}
By definition $\Delta$ is linear in the second argument and by Definition \ref{def:contra-gene} (ii)  is symmetric. Hence $\Delta$ is bilinear and by
 Definition \ref{def:contra-gene} (i)  $\Delta(P,P)=0$ for all
$P\in V_F$. This proves (i) and (ii) is proved similarly. 
\end{demo}

\begin{corollary}
Let $n\in V_\F^\ast$. Then
$\F_n$ is an affine space for ${\cal S}_0^2(V_\F^\ast)$. In particular \rm{Card}$(\F_0)=8$.
\end{corollary}

\begin{demo}
If $n=0$ this is exactly the proposition above.
 For arbitrary $n \in V_F^\ast$ the statement follows
from the proposition above and Proposition \ref{prop:corres}.
\end{demo}

\section{Oriented maps and octonion multiplication}
In this section  by ``exponentiation'' we establish 
 an equivalence between oriented (resp. $n-$oriented, $n\ne0$) maps on the Fano plane $\F$
and octonion (resp. split octonion) multiplication laws on ${\mathbb O_\F}$.  We also prove the converse to Theorem \ref{theo:future} (see Theorems \ref{theo:future1} and \ref{theo:future2}).

\begin{definition}\label{def:exp}
Let $\F$ be a Fano plane and let $(P\mapsto \alpha_P)$ be an oriented map. Define the associated multiplication factors: 
$\pm e^{\alpha}:\F_0^2\to \{-1,1\}$  by
\beqa
\nn
\pm e^{\alpha}_{PQ} = e^{\alpha_P(Q)}\quad\forall (P,Q)\in \F_0^2.
\eeqa

\end{definition}

\begin{proposition}\label{theo:exp}
 Let $\F$ be a Fano plane, 
 let $(P\mapsto \alpha_P)$ be an oriented map and let $\pm  e^{\alpha}$  be the associated multiplication factors. Then:
\begin{enumerate}
\item[(i)]  $({\mathbb O_\F}, {\bf 1},  e^{\alpha}) $ is  a composition algebra
and   $({\mathbb O_\F}, {\bf 1},  e^{\alpha}) \in {\mathbb O_\F}^+_{1}$.
\item[(ii)]  $({\mathbb O_\F}, {\bf 1}, - e^{\alpha})$ is a composition algebra
and   $({\mathbb O_\F}, {\bf 1},  -e^{\alpha}) \in {\mathbb O_\F}^-_{1}$.
\end{enumerate}
\end{proposition}

\begin{demo}
To  prove (i) we have to show (see Proposition \ref{prop:comp}) that
\beqa
(1)&:& e^{\alpha}_{PQ} e^{\alpha}_{Q P+Q}=1\ , \  \  \forall P\ne Q \in \F\ , \nn\\
(2)&:& e^{\alpha}_{PQ} e^{\alpha}_{QR} e^{\alpha}_{RS} e^{\alpha}_{SP} = -1 \ , \ \  \mbox{if} \ \ \{P,Q,R,S\} \ \ \mbox{is a quadrilateral}\nn
\eeqa
which is equivalent to
\beqa
e^{\alpha_P(Q)} e^{\alpha_Q(P+Q)}&=&1 \ , \nn\\
e^{\alpha_P(Q)} e^{\alpha_Q(R)} e^{\alpha_R(S)} e^{\alpha_S(P)} &=&-1 \ , \nn
\eeqa
and, taking logarithms, this is equivalent to
\beqa
\label{eq:alpha}
\alpha_P(Q) + \alpha_Q(P+Q)&=&0\ , \nn\\
\alpha_P(Q) + \alpha_Q(R) + \alpha_R(S) + \alpha_S(P) &=&1
 \ . 
 \eeqa
 The first equation reduces to
 \beqa
\alpha_P(Q) + \alpha_Q(P)+\alpha_Q(Q)=0 
\eeqa
and this follows from Definition \ref{def:contra-gene} (i) and (ii).

To prove the second equation recall that $P+Q+R+S=0$ if $\{P,Q,R,S\}$ is a quadrilateral ({\it c.f.} L7 in List \ref{li:inci}). Hence
\beqa
\alpha_P(Q) + \alpha_Q(R) + \alpha_R(S) + \alpha_S(P)&=& \phantom{+}\alpha_P(Q) +  \alpha_Q(R) + \alpha_R(P+Q+R) + \alpha_{P+Q+R}(P)\nn\\
&=& \phantom{+}\alpha_P(Q) + \alpha_Q(R)\nn\\
&&+ \alpha_R(P)+ \alpha_R(Q)+ \alpha_R(R)\nn\\
&&+1+\alpha_P(P) + \alpha_P(Q) +\alpha_P(R) =1\nn
\eeqa
by  Definition \ref{def:contra-gene} (i) and (ii).

Finally  if $P\in \F$ then $\overrightarrow{P}$ is a line because it is the set of zeros of the linear form $\alpha_P$ and
this proves that $({\mathbb O_\F}, N_{\mathbb O_\F},e^{\alpha}) \in {\mathbb O_\F}_1^+$.

To prove (ii) it is sufficient to remark that the affine functions  $\alpha'_P$
\beqa
\alpha'_P(Q)= \alpha_P(Q) +1 \hskip .7 truecm  \forall (P,Q) \in \F_0^2 \ \nn
\eeqa
satisfy \eqref{eq:alpha} and that
$\overrightarrow{P_{e^{\alpha}}}= \overleftarrow{P_{-e^{\alpha}}}$.

\end{demo}

We now show that  all multiplication factors which   define a composition law on ${\mathbb O_\F}$ are uniquely obtained in one of  these ways.
The crucial point  here is Theorem  \ref{theo:future}.
 
\begin{proposition}\label{prop:log} 
 Let $({\mathbb O_\F}, {\bf 1}, \epsilon)$ be a composition algebra.
\begin{enumerate}
\item[(i)] If   $({\mathbb O_\F}, {\bf 1}, \epsilon) \in {\mathbb O_\F}^+_{1}$ then there exists a unique oriented map
$(P \mapsto \alpha_P)$ such that $\epsilon=e^{\alpha}$. 

\item[(ii)]If   $({\mathbb O_\F}, {\bf 1}, \epsilon) \in {\mathbb O_\F}^-_{1}$ then there exists a unique  oriented map
$(P \mapsto \alpha_P)$ such that  $\epsilon=-e^{\alpha}$. 
\end{enumerate}
\end{proposition}

\begin{demo}
(i):  For all $(P,Q)\in\F^2$ define 
\beqa
\alpha_P(Q)= \left\{
\begin{array}{cll}
\ell (\epsilon_{PQ})&\mbox{if}& Q\ne P,\\
1&\mbox{if}& Q= P \ . 
\end{array}
\right.\nn
\eeqa
To prove (i) we have to show that $(P \mapsto \alpha_P)$ is an oriented map. If $P\ne Q$, since $\epsilon_{PQ}= -\epsilon_{QP}$ it is clear that
$\alpha_P(Q) + \alpha_Q(P)=1$.
It remains  to show that $\alpha_P$ is a linear function of $Q$ and for this
it is sufficient to show that $\alpha_P$ has exactly {\it three zeros in $\F$ which are aligned}.
 Note that  if $Q \not =P$ then: $\alpha_P(Q)=0$ {\it iff} $\epsilon_{PQ}=1$. Since $\alpha_P(P)=1$ it follows that
 \beqa
\{Q \in \F  \ \ \text{s.t.} \ \ \alpha_P(Q)= 0 \} = \{ Q\in\F \ \ \text{s.t.} \ \ \epsilon_{PQ}=1\} =
\overrightarrow{P} \ . \nn
\eeqa
However  $\overrightarrow{P}$ is a line in $\F$ by Theorem \ref{theo:future} and so  $\alpha_P$ is  linear .
The proof of (ii) is analogous.
\end{demo}

The following two theorems are converses to the two parts of Theorem \ref{theo:future}.
{
\begin{theorem} \label{theo:future1}
Let $\F$ be a Fano plane and let  $\epsilon$ be a multiplication factor. Then $({\mathbb O_\F}, {\bf 1}, \epsilon)$ 
is a composition algebra if either  $\overrightarrow{P}$   is a line for all $P$ in $\F$ or  $\overleftarrow{P}$  is a line for all $P$ in $\F$.
\end{theorem}
\begin{demo} Suppose $\overrightarrow{P}$   is a line for all $P$ in $\F$.
 The proof of Proposition \ref{prop:log}(i) above shows that there exists an oriented map $(P \mapsto \alpha_P)$ such that
 $\forall P\ne Q \in \F$,
\beqa
\epsilon_{PQ}=e^{\alpha_P(Q)}
\eeqa
since the only property (besides antisymmetry) of $\epsilon$  used in the proof was that $\overrightarrow{P}$   is a line for all $P$  in $\F$. It then
follows from Proposition \ref{theo:exp}(i) that $({\mathbb O_\F}, {\bf 1}, \epsilon)$ is a composition algebra.

\end{demo}
\begin{theorem} \label{theo:future2}
Let $\F$ be a Fano plane and let  $\Delta$ be a set of seven triangles satisfying the axioms
of a Fano plane. Then there  exists a unique multiplication factor $\epsilon$ such that $({\mathbb O_\F}, {\bf 1}, \epsilon)$ 
is a composition algebra and $\Delta=\{\overleftarrow{P} \ \ \text{s.t.} \ \ P\in\F\}$.

\end{theorem}
}
\begin{demo}
By Proposition \ref{prop:trian} there exists a unique oriented map $P\mapsto \alpha_P$ such that
\beqa
\Delta= \Big\{\Delta_P \ \ \text{s.t.} \ \ P\in \F\Big\} \nn
\eeqa
where $\Delta_P =\F\setminus\{P\}\setminus \big\{Q \in \F \ \ \text{s.t.} \ \ \alpha_P(Q)=0\big\}$.
Define the multiplication factor $\epsilon$ by $\epsilon= e^\alpha$. Then $(\mathbb O, {\bf 1}, \epsilon)$ is a composition algebra by Proposition \ref{theo:exp} and by definition
\beqa
\overleftarrow{P}= \Big\{Q \in \F \ \ \text{s.t.} \ \ \epsilon_{PQ} =-1\Big\}
= \Big\{Q \in \F\setminus\{P\} \ \ \text{s.t.} \ \ \alpha_P(Q)=1 \Big\}= \Delta_P \ . \nn
\eeqa
\end{demo}

We now show how the above results for division composition algebras $(\mathbb O_\F, {\bf 1}, \epsilon)$
can be used to deduce analogous result for split composition algebras $(\mathbb O_\F, e^{n'}, \epsilon')$.

\begin{lemma}\label{lem:comp-spli}
\mbox{\color{white} t}
\begin{enumerate}
\item[(a)]
Let  $n' \in V_\F^\ast$ and let $({\mathbb O_\F}, e^{n'}, \epsilon')$ be a composition algebra. For all  $P\neq Q \in \F$ define
\beqa
\epsilon_{PQ}= e^{n'(P)n'(Q)} \epsilon'_{PQ} \ .\nn
\eeqa
Then  $({\mathbb O_\F}, {\bf 1},\epsilon )$ is a composition algebra.
\item[(b)] Let $({\mathbb O_\F}, {\bf 1}, \epsilon)$ be a composition algebra.
Let $n' \in V_\F^\ast$ and
for all  $P\neq Q \in \F$ define
\beqa
\epsilon'_{PQ}= e^{n'(P)n'(Q)} \epsilon_{PQ} \ .\nn
\eeqa
Then  $({\mathbb O_\F}, e^{n'},\epsilon )$ is a composition algebra.
\end{enumerate}
\end{lemma}

\begin{demo}
We only prove (a) since the proof of $(b)$ is the same.
By \eqref{eq:comp1} and \eqref{eq:comp2} we have to show that
\begin{enumerate}
\item[(i)] $\epsilon_{PQ} \epsilon_{Q(P+Q)} =1$ for all $P\neq Q \in \F$;
\item[(ii)]  $\epsilon_{PQ} \epsilon_{QR}  \epsilon_{RS} \epsilon_{SP} =-1$ for every quadrilateral $\{P, Q,R,S\}$.
\end{enumerate}
(i): We have
\beqa
\epsilon_{PQ} \epsilon_{Q(P+Q)} = e^{n'(P)n'(Q) + n'(Q)n'(P+Q)} \epsilon'_{PQ} \epsilon'_{Q(P+Q)}= e^{n'(Q)} \epsilon'_{PQ} \epsilon'_{Q(P+Q)}\ .
\eeqa
However by \eqref{eq:comp1} we have 
\beqa
\epsilon'_{PQ} \epsilon'_{Q(P+Q)}= e^{n'(P+P+Q)} = e^{n'(Q)} \  \nn
\eeqa
which proves (i).

\noi
(ii):  We have
\beqa
\epsilon_{PQ} \epsilon_{QR}  \epsilon_{RS} \epsilon_{SP} = e^{n'(P)n'(Q) + n'(Q)n'(R)+n'(R)n'(S) + n'(S) n'(P)} \epsilon'_{PQ} \epsilon'_{QR}  \epsilon'_{RS} \epsilon'_{SP} \ . \nn
\eeqa
Using
\beqa
\epsilon'_{PQ} \epsilon'_{QR}  \epsilon'_{RS} \epsilon'_{SP}= - e^{n'(P+Q) + n'(P+S)} \nn
\eeqa
and $S=P+Q+R$ this simplifies to $-1$ which proves (ii) and this completes the proof of (a).
      \end{demo}

\begin{corollary}
Let $({\mathbb O_\F}, {\bf 1}, \epsilon)$ be a composition algebra.  Let $n' \in (V_{\F})^\ast \setminus \{0\}$ and let 
$$
\epsilon'_{PQ}= e^{n'(P)n'(Q)} \epsilon_{PQ} \ \ \ , \ \ \forall P\ne Q \in \F
$$
so that  (Lemma \ref{lem:comp-spli}) $({\mathbb O_\F}, e^{n'}, \epsilon')$ is a composition algebra.
\begin{enumerate}
\item[(i)] If $({\mathbb O_\F}, {\bf 1}, \epsilon) \in {\mathbb O_\F}_1^+$  define   $P_{n'} \in \F$ by
$\overrightarrow{P_{n'}}=\mbox{Ker}(n')$. Then  $P_{n'}$ is the unique point of $\F$ such that
$\overrightarrow{P_{n'}} =\F \setminus \{P_{n'}\}$.
\item[(ii)]If   $({\mathbb O_\F}, {\bf 1}, \epsilon) \in {\mathbb O_\F}_1^-$  define   $P_{n'} \in \F$ by
$\overleftarrow{P_{n'}}=\mbox{Ker}(n')$.  Then  $P_{n'}$ is the unique point of $\F$ such that
$\overleftarrow{P_{n'}} =\F \setminus \{P_{n'}\}$.
\end{enumerate}
\end{corollary}
\begin{demo}
Let $P \in \F$. Then we have the following sequence of equivalences:
\beqa
&&\epsilon'_{PQ} = 1 \ \ \forall Q \ne P\nn\\
\Leftrightarrow&& e^{n'(P)n'(Q)} \epsilon_{PQ}=1 \ \ \forall Q \ne P\nn\\
\Leftrightarrow&&\left\{\begin{array}{l}
                   n'(P) n'(Q)  = 0 \ \ \forall Q \in \overrightarrow{P}\ , \\
                    n'(P) n'(Q) = 1 \ \ \forall Q \in \overleftarrow{P}
                    \end{array}\right.\nn\\
 \Leftrightarrow&&\left\{\begin{array}{l}
                    n'(P)=1\ , \nn\\
                     n'(Q)  = 0 \ \ \forall Q \in \overrightarrow{P}\ , \\
                      n'(Q) = 1 \ \ \forall Q \in \overleftarrow{P}
                      \end{array}\right.\nn\\
\Leftrightarrow&&  \overrightarrow{P}= \mbox{Ker}(n') \ . \nn                     
\eeqa
This proves (i) and (ii) is proved similarly.
\end{demo}
\begin{definition}\label{def:fin}
Let $n'\in V_\F^\ast\setminus\{0\}$. We denote by ${\mathbb O_\F}_{e^{n'}}^+$  (resp. ${\mathbb O_\F}_{e^{n'}}^-$) the set of all composition algebras $({\mathbb O_\F}, e^{n'}, \epsilon')$ such that there is a unique point
$P_{n'} \in \F $  satisfying  $\overrightarrow{P_{n'}}=\F\setminus\{P_{n'}\}$ (resp. $\overleftarrow{P_{n'}}=\F\setminus\{P_{n'}\}$).
\end{definition}
We conclude this paper with the main theorem which shows the equivalence between  $n-$oriented Fano planes $(\F,n,\alpha)$  (see Definition \ref{def:contra-gene}),
  and composition algebras $({\mathbb O_\F}, N,\e)$ (see Definition \ref{def:mo}) such that the restriction of $N$ to
$\F\subset {\mathbb O_\F}$ is a norm on $\F$ (see Definition \ref{def:norm}).

\begin{theorem} \label{theo:div-split}
Let $\F$ be a Fano plane and let $n' \in V_ \F^\ast$.
\begin{enumerate}
\item[(i)] If  $\alpha \in \F_{n'}$   then $({\mathbb O_\F}, e^{n'}, e^{\alpha}) \in {\mathbb O_\F}_{e^{n'}}^+$ and  $({\mathbb O_\F}, e^{n'}, -e^{\alpha}) \in {\mathbb O_\F}_{e^{n'}}^-$.
\item[(ii)] If $({\mathbb O_\F}, e^{n'},\e) \in {\mathbb O_\F}_{e^{n'}}^+$ (resp.   $({\mathbb O_\F}, e^{n'},\e) \in
{\mathbb O_\F}_{e^{n'}}^-$) then  there  exists a unique    $\alpha \in \F_{n'}$  such that 
 $\epsilon= e^\alpha$ (resp. $\epsilon= -e^{\alpha}$).

\end{enumerate}

\end{theorem}

\begin{demo}
In the case $n'=0$, part (i) is Proposition \ref{theo:exp} and part (ii) is Proposition \ref{prop:log}.
If $n'\in V_\F^\ast$ recall that $\F_{n'}$ is the set of all $n'-$oriented maps (see Definition
\ref{def:contra-gene}) and for  ${\mathbb O_\F}^\pm_{e^{n'}}$ see Definition \ref{def:Opm} ($n'=0$) and
Definition \ref{def:fin} ($n'\ne 0$).
We only prove the theorem for $\F_{n'} \slash \mathbb {O_\F}_{e^{n'}}^+$, the proof for   $\F_{n'} \slash{\mathbb O_\F}_{e^{n'}}^-$ being entirely analogous. 
Consider the diagram
\beqa
\xymatrix{ \F_0 \ar[rr]^E \ar[dd]_A  && {\mathbb O_\F}^+_{1}\ar[dd]^B\\ \\
      \F_{n'} \ar@{.>}[rr]_{E' } &&{\mathbb O_\F}^+_{e^{n'}}
      }
\eeqa
where
\begin{enumerate}
\item  $A$ is the bijection of Proposition \ref{prop:corres};
\item $E$ is the bijection of Proposition  \ref{prop:log} (i);
\item $B$ is the bijection Lemma \ref{lem:comp-spli} (b);
\item $E'= B \circ E \circ A^{-1}$.
\end{enumerate}

\noi
It is clear that $E'$ is a bijection
so it only remains to show that if $(\alpha',n') \in \F_{n'}$ then $E'(\alpha',n') = (\mathbb O_\F, e^{n'},\epsilon'')$  where
\beqa
\epsilon''_{PQ} = e^{\alpha'_P(Q)} \ , \  \ \forall P\ne  Q \in \F  . \nn
\eeqa

Let $(\alpha',n')\in \F_{n'}$. Then $A^{-1}(\alpha',n')= \alpha$ where
\beqa
\alpha_P(Q) =\alpha'_p(Q)+ n'(P) n'(Q) \ , \  \ \forall P\ne  Q \in \F  \ .\nn
\eeqa
Next  we have $E(\alpha) =(\mathbb O_\F, {\bf 1}, \epsilon)$ where
\beqa
\epsilon_{PQ}= e^{\alpha'_P(Q) + n'(P) n'(Q)} \ , \   \ \forall P\ne  Q \in \F  \ ,\nn
\eeqa
and finally we have $B(\mathbb O_\F, {\bf 1}, \epsilon)= (\mathbb O_\F, e^n, \epsilon')$ where
\beqa
\epsilon'_{PQ} = e^{n'(P) n'(Q)} \epsilon_{PQ}   \ , \   \ \forall P\ne  Q \in \F  \ .\nn
\eeqa
It is now clear that $\epsilon'=\epsilon''$.
\end{demo}

\bibliographystyle{utphys}
\bibliography{ref}

\end{document}